\documentclass[12pt,a4paper]{article}
\usepackage[cp1251]{inputenc}
\usepackage{amsmath}
\usepackage{amssymb}
\usepackage{latexsym}
\usepackage{srcltx}
\usepackage{graphics}
\usepackage{color}
\usepackage{epsfig}

\newtheorem{theorem}{Theorem}
\newtheorem{prop}{Proposition}
\newtheorem{lemma}{Lemma}
\newtheorem{cor}{Corollary}
\newtheorem{remark}{Remark}

\newcommand{\re}{{\mathbb R}}

\newcommand{\T}{{\mathbb T}}

\newcommand{\cA}{{\cal{A}}}

\newcommand{\cB}{{\cal{B}}}
\newcommand{\cH}{{\cal{H}}}

\newcommand{\Tr}{\mathrm{Tr}}

\date{}

\author{V.Yu.Protasov\thanks{Lomonosov Moscow
State University, National Research University High School of Economics, {e-mail: \tt\small
v-protassov@yandex.ru}}, M.E.Shirokov\thanks{Steklov Mathematical Institute, {e-mail: \tt\small
msh@mi-ras.ru}}}

\title{On mutually inverse transforms
of functions on a half-line}

\begin{document}
\maketitle

\begin{abstract}

Two transforms of functions on a half-line are considered. It is proved that their composition gives
a concave majorant for every nonnegative function. In particular, this composition is the identity transform on
the class of nonnegative concave functions. Applications of this result to some problems of mathematical physics are indicated.
Several open questions  are formulated.
\smallskip
\end{abstract}
\bigskip

\begin{center}
\textbf{The main result}
\end{center}

\medskip

The following transforms of functions on the positive
half-line $\mathbb{R}$ naturally arise in problems of mathematical
physics (see Section 3):
$$
F[f](x)\ = \ \sup_{\, t > 0}\ \frac{f(xt)}{t+1}\ ,
\qquad G[f](x)\ = \ \inf_{\, t > 0}\  f(xt)\, \Bigl(1 \, + \, \frac1t  \Bigr).
$$
Throughout this paper, we assume that these transforms
are applied to functions on which they are well
defined, i.e., $F[f](x)$ and  $G[f](x)$ are finite for any $x$.
For example, for $F[f]$, it suffices that the function $f(t)$
be at most linearly growing, i.e.,
$$
\sup_{t \in [0,a]} f(t) \, < \, C(a+1), \, a > 0.
$$
The composition is denoted by $G[F[f]] = GF[f]$.
For concave functions $f$ of a certain class, it was shown
in \cite{ECN} that the composition is an identity transform.
We will prove that $GF[f] = f$ if and only if $f$ is any
monotone concave function satisfying some conditions. Moreover, for an arbitrary nonnegative
function $f$, the composition $GF[f]$ is the
smallest concave majorant of $f$, i.e., it coincides with the
function
$\ \tilde f(t) \ = \ \inf\, \Bigl\{
\varphi(t)\ \Bigl| \ \varphi \ge f\, , \ \varphi \ \mbox{concave}\, \Bigr\}.
$
This property of $F$ and $G$ resembles the Legendre–Young duality transform~\cite{MT}, but no relations to it have
been found.

Consider two classes of concave functions on $\re_+$.
The class $\cA$ consists of concave nondecreasing nonnegative
functions. The graph of each function $\varphi \in \cA$
consists of three parts (from left to right): the segment
from the point $(0,0)$ to some point $A_1$, the arc of the
concave function from $A_1$ to $A_2$, and the horizontal ray
from $A_2$ to $+\infty$. Some of these parts may be degenerate,
which corresponds to the cases $A_1 = (0,\varphi(0))$, $A_2 = A_1$ or $A_2 = +\infty$.
For example, the function $f\equiv1$ belongs to $\cA$ and for it $A_1 = A_2 = (0,1)$.
Any concave increasing
function satisfying the condition $\varphi(0) > 0$
belongs to $\cA$, and, for it $A_1 = (0,\varphi(0)), A_2 = +\infty$. The
class $\cB$ consists of concave nonincreasing nonpositive
functions with the asymptote $y = kt, \ t\to +\infty$. The
graph of such a function consists of three parts (from
left to right): a horizontal segment to some point $B_1$,
the arc of the concave function from $B_1$ to $B_2$, and the
ray from $B_1$ to $+\infty$, whose extension beyond the point $\cB_2$
passes through zero. Some of these parts may be
degenerate. For example, for $f\equiv -1$ we have
$B_1 = B_2 = +\infty$; while for the function  $\varphi(t) = - \frac{1}{t+1} - t$ we have
$B_1 = (0,-1), B_2 = +\infty$.

The classes $\cA$ and $\cB$ are convex and closed in the
topology of uniform convergence on each compact set.
They intersect only on the zero function. The $\cA$-hull
of a function $f$ is the smallest majorant by functions from $\cA$, i.e.
$f_{\cA}(t) = \inf\, \bigl\{\varphi(t) \ \bigl| \  \varphi \ge f\, , \ \varphi \in \cA \bigr\}$. Similarly we define the $\cB$-hull: $f_{\cB}(t) = \inf\, \bigl\{\varphi(t) \ \bigl| \  \varphi \ge f\, , \  \varphi \in \cB \bigr\}$. Note that the $\cA$-hull is defined for any function $f$ with at
most linear growth, while the $\cB$-hull is defined only
for non-positive $f$.

\begin{theorem}\label{th.10}
If a function  $f$ takes at least one positive value then ${GF[f] = f_{\cA}}$, otherwise $GF[f] = f_{\cB}$.
\end{theorem}
\begin{cor}\label{cor.10}
$GF[f] = f$ if and only if  $f \in \cA\cup \cB$.
\end{cor}
\begin{cor}\label{cor.20}
If $f$ is a nonnegative function then  $GF[f]$ is the concave majorant of $f$. For nonnegative functions,
the equality $GF[f] = f$ holds if and only if $f$ is
concave.
\end{cor}
The proof of Theorem 1 is geometric. We begin
with an auxiliary assertion. Suppose that the function
$f: \re_+ \to \re$ grows at most linearly. Through a given point $A = (-a,0), \, a> 0$
we draw a minimum-slope straight line $\ell$
lying above the graph of f, i.e. $\ell(t) \ge f(t), \ t \in \re_+$
(here and below, a linear function is denoted the same
as the corresponding line). We call a left supporting
line of the graph of $f$. In contrast to a usual \emph{supporting
line}~$\ell$, for which $\inf_{t\in \re_+}\bigl(\ell(t) - f(t)\bigr) = 0$, the left supporting
line must intersect the negative horizontal half-line ${\{(-a, 0)\, | \, a\ge 0\}}$.
\begin{lemma}\label{l.10}
For any $x > 0$, the straight line passing
through the points $(-x, 0)$ and $(0, F[f](x))$ is a left supporting
line to the graph of the function $f$.
\end{lemma}
{\tt Proof.} At the point with abscissa $z$, this line has the
ordinate
$$
y = \frac{(x+z) F[f](x)}{x} = (1+t)F[f](x),
$$
where $t = z/x$. Hence,
$$
\inf_{z \ge 0}
\bigl( y(z) - f(z)\bigr) \, = \, \inf_{t \ge 0}
\bigl( (1+t)F[f](x) - f(xt)\bigr) \, \ge \, 0.
$$
Therefore, the line is above the graph of $f$ and its slope
cannot be reduced, since $\inf_{t \ge 0}
\bigl( F[f](x) - \frac{1}{1+t}\, f(xt)\bigr) \, = \, 0$.
{\hfill $\Box$}
\smallskip

\noindent {\tt Proof of the theorem~\ref{th.10}}. Denote  $z  = xt$. Then
$G[f](x)\ = \ \inf_{z>0}\, f(z)\, \frac{x+ z}{z}\, . $
Hence,
\begin{equation}\label{eq.G1}
GF[f](x)\ = \ \inf_{z>0}\, F[f](z)\, \frac{x+ z}{z}\, .
\end{equation}
From the point $A=(-z,0)$, we draw a supporting line
to the graph of $f$. It intersects the vertical axis at the
point $M = (0, F[f](z))$. Let also $B = (x, 0)$. Drop a
perpendicular to the horizontal axis at the point $B$; it
crosses the supporting line at the point $N$. Since the
triangles $AMO$ and $ANB$, where $O$ is the origin, are similar, we obtain.
$BN/OM = AB/AO = (z+x)/z$. Thus,
$BN = F[f](z)(z+x)/z$. Combining this with~(\ref{eq.G1}), we conclude that $GF[f](x)$ is equal
to the smallest value of $\ell(x)$ over all left supporting lines $\ell$ to the
graph of $f$. Let denote $\tilde f$ the smallest concave majorant
of $f$. It is the pointwise minimum of all supporting lines
to the graph of $f$.

Assume that $f(t) > 0$ at least at one
point $t$. In the graph of $\tilde f$ there are points $A_1, A_2$
(which may coincide or grow to $+\infty$) such that the supporting
line to the graph of $\tilde f$ at the point $A_1$ passes
through the origin (we use the rightmost of these
points), while the supporting line at the point $A_2$ is
horizontal (we use the leftmost point). Thus, $A_2$ is a
maximizer of $\tilde f$ and, by assumption, this maximum is
positive. For any point between $A_1$ and $A_2$, the supporting
line to $\tilde f$ intersects the negative horizontal
half-line; hence, this is a left supporting line. Therefore,
the graph of $GF[f]$ coincides with the graph of $\tilde f$
between the points $A_1$ and $A_2$, with the graph of the left
supporting line through the origin to the left of $A_1$, and
with the horizontal line through to the right of $A_2$.
Thus, $GF[f] = f_{\cA}$. The second case is treated in a
similar manner: $f \le 0$ then $GF[f] = f_{\cB}$.
{\hfill $\Box$}
\smallskip

\begin{center}
\textbf{2. Comments}
\end{center}

\begin{remark}\label{r.10}{\em After interchanging the factors in the
composition, the resulting operator $FG$ is not an identity
one in the classes $\cA$ and $\cB$. We do not know, for
which functions $f$, it holds that $FG[f] = f$. One condition
can be obtained by applying the following argument.
Let denote $\ell_z$ the line passing through the
points $(-z, 0)$ and $(0, f(z))$.  It is easy to see that
$G[f](t) = \inf_{z> 0} \ell_z(t)$. Then the equality $FG[f] = f$ is
equivalent to the fact that each line of the family
$\{\ell_z \, | \, z > 0\}$  is a right supporting line to the graph of $G[f]$. If $f$
is nonnegative, this condition is equivalent to this
requirement for $f$. Denote by $\cH$ the family of functions $\varphi(t) = \frac{at}{t+b}\,  , \ a, b > 0$.
The equality $FG[f] = f$ holds if and only if $G[f]$ is not identically zero and, for any pair
of points in the graph of $f$, the arc of a function from
the family $\cH$ that connects these points lies above the
graph of $f$.
Finding more adequate conditions for the equality $FG[f] = f$
remains an open problem.}
\end{remark}

\begin{center}
\textbf{3. Some applications}
\end{center}

Denote by $\T$ the convex cone in a linear space. Let
$\alpha,\beta$ and $\gamma$  be nonnegative linear functions on $\T$
such that $\alpha(t)>0$ for all $t$ from  $\T$ not equal to zero.
Consider the nonnegative functions
$$
f(x)=\sup\{\gamma(t)\,|\,\alpha(t)\leq 1,\beta(t)\leq x\},\quad g(x)=\sup\{\gamma(t)\,|\,\alpha(t)+x^{-1}\beta(t)\leq 1\}
$$
on $\mathbb{R}_+$. Corollary 1 allows to prove the following

\begin{prop}\label{p1} The above functions are connected by
the relations $\,g=F[f]$ and $f=G[g]$, where $F$ and $G$
are the transforms defined in Section 1.
\end{prop}

{\tt Proof.} It is easy to see that the first function is concave.
Therefore, by Corollary 1, it suffices to show that
\begin{equation}\label{t}
g(x)=\sup_{r\in(0,1)} rf\left(x\,\frac{1-r}{r}\right)=\sup_{t>0} \frac{f(xt)}{1+t}.
\end{equation}

To prove (\ref{t}) note first that $g(x)\geq rf\left(x\frac{1-r}{r}\right)$ for all $r\in(0,1]$. For any $\varepsilon>0$ there is $t_{\varepsilon}\in \T$ such that  $\alpha(t_{\varepsilon})\leq 1,\beta(t_{\varepsilon})\leq x\frac{1-r}{r}$ and $f\left(x\frac{1-r}{r}\right)\leq\gamma(t_{\varepsilon})+\varepsilon$. Then
$rt_{\varepsilon}\in \T$ and $\alpha(rt_{\varepsilon})+x^{-1}\beta(rt_{\varepsilon})\leq r+(1-r)\leq 1$. Hence $g(x)\geq \gamma(rt_{\varepsilon})\geq rf\left(x\frac{1-r}{r}\right)-r\varepsilon$. Since $\varepsilon$ is arbitrary, we obtain the required inequality.

For any $\varepsilon>0$ there is $t_{\varepsilon}\in \T\setminus\{0\}$ such that  $\alpha(t_{\varepsilon})+x^{-1}\beta(t_{\varepsilon})\leq 1$ and $g(x)\leq \gamma(t_{\varepsilon})+\varepsilon$. Since  $r=\alpha(t_{\varepsilon})\in(0,1]$,   $\alpha(t_{\varepsilon}/r)=1$.
Then  $\beta(t_{\varepsilon}/r)\leq x(1-r)/r$ and hence $\gamma(t_{\varepsilon}/r)\leq f\left(x\frac{1-r}{r}\right)$. Thus,
$$
g(x)\leq \gamma(t_{\varepsilon})+\varepsilon\leq rf\left(x\frac{1-r}{r}\right)+r\varepsilon.
$$
It suffices to note that the condition $r\in(0,1)$ in (\ref{t}) can be replaced by the condition $r\in(0,1]$ by the concavity of $f$.
{\hfill $\Box$}\smallskip

\textbf{Example 1.}  Let $\T$ be the cone of positive operators
on a Hilbert space $H$ with the inner product $\langle\cdot,\cdot\rangle$ and the norm $|\cdot|$.
To avoid the technical difficulties, we assume that
$H$ is finite-dimensional.

Given any positive operator $G$ on $H$ with a zero
lower spectral bound and any $E>0$, the space $H$ can
be equipped with the inner product
$$
\langle\varphi,\psi\rangle^G_E=\langle\varphi,\psi\rangle+\langle\varphi,G\psi\rangle/E,\quad \varphi,\psi\in H
$$
and the corresponding norm $|\varphi|^G_E=\sqrt{\langle\varphi,\varphi\rangle^G_E}$ of any vector $\varphi\in H$. Then, for any linear operator $A$ on $H$, the
quantity
$$
|||A|||^G_E=\sup\{ |A\varphi|\,|\, \varphi\in H,\, |\varphi|^G_E\leq 1 \}
$$
is the operator norm of $A$ regarded as an operator from
the space $H$ with the inner product $\langle\cdot,\cdot\rangle^G_E$ to the space $H$ with the inner product $\langle\cdot,\cdot\rangle$. One can show (see \cite{ECN}) that
$$
|||A|||^G_E=\sup\{\sqrt{\Tr AXA^*} \,|\, X\in \T,\, \Tr X+E^{-1}\Tr GX  \leq 1\}.
$$
For any linear operator $A$ on $H$, we can define the family
of norms
$$
||A||^G_E=\sup\{ |A\varphi|\,|\, \varphi\in H,\, |\varphi|\leq 1,\; \langle\varphi,G\varphi\rangle\leq E \},\quad E>0,
$$
which are useful for a number of problems of mathematical
physics \cite{ECN}.\footnote{Specific properties of these norms are exhibited in the case of an
unbounded operator $G$ and an infinite-dimensional Hilbert
space $H$, in which these norms define topologies weaker than
the topology of the operator norm.} It is proved in \cite{ECN} that
$$
\|A\|^G_E=\sup\{\sqrt{\Tr AXA^*} \,|\, X\in \T,\, \Tr X\leq 1, \Tr GX  \leq E\}.
$$
This representation shows that the function
$E\rightarrow[\|A\|^G_E]^2$ is concave on $\mathbb{R}_+$ and tends to the ordinary operator norm of $A$ as $E\rightarrow+\infty$.
By Proposition 1 for each operator $A$ on $H$ the nonnegative functions $f_A(E)=\left[\|A\|^G_E\right]^2$ and
$g_A(E)=\left[|||A|||^G_E\right]^2$ are related by the transforms $F$ and $G$:
\begin{equation}\label{u-rel}
g_A=F[f_A],\quad \quad f_A=G[g_A].
\end{equation}
These relations were obtained in \cite[Theorem 3]{ECN} by applying a
rather complicated method involving special properties
of the functions $f_A$ and $g_A$. Proposition 1 shows that
these relations are a special case of the more general
property presented in Corollary 1.

The norms $\|\cdot\|^G_E$  (naturally defined for operators on
separable Hilbert spaces in the case of an unbounded
operator $G$) were found useful for some problems
in the theory of open quantum systems. They
were used to prove a generalization of the Kretschmann–
Schlingemann–Werner theorem \cite[Theorem 2]{ECN}, which made it possible to obtain a Stinespring
representation for strongly converging
sequences of quantum channels. These norms also
arise in describing strongly continuous quantum
dynamic semigroups in the case of energy constraints
imposed on the states of the quantum system \cite{QDS}.

The norms $|||\cdot|||^G_E$ naturally arise in the theory of relatively
bounded operators on infinite-dimensional
Hilbert spaces \cite{Kato, BS}.

Relation (\ref{u-rel})  defines a one-to-one correspondence
between the norms $\|\cdot\|^G_E$ and $|||\cdot|||^G_E$ as functions on $\mathbb{R}_+$,
which is of interest for both theoretical considerations
and physical applications (see \cite{ECN}, the example after
Theorem 3).


\bigskip


\begin{thebibliography}{NN}



\bibitem{MT}
G.G.Magaril-Il’yaev, V.M.Tikhomirov,
"Convex
Analysis: Theory and Applications",  Am. Math. Soc. Providence, R.I., 2003.

\bibitem{ECN} M. E. Shirokov, "Operator E-norms and their use", Sb. Math., 211:9 (2020), 1323–1353; arXiv:1806.05668.

\bibitem {QDS} A.S.Holevo, M.E.Shirokov, "Energy-constrained diamond norms and quantum dynamical semigroups", Lobachevskii J. Math., 40:10 (2019), 1569–1586; arXiv:1812.07447.

\bibitem{Kato} T.Kato, "Perturbation Theory for Linear Operators". Springer-Verlag, New York-Heidelberg-Berlin, 1980.

\bibitem{BS}  B.Simon, "Operator Theory: A Comprehensive Course in Analysis, Part IV",  American Mathematical Society, 2015.
\end{thebibliography}
\end{document}